\documentclass[11pt,leqno]{article}
\usepackage[english]{babel}
\usepackage{tikz}
\usepackage{graphicx}
\usepackage{booktabs,tabularx}
\usepackage{verbatim}
\usepackage{etoolbox}
\usepackage{tkz-euclide}
\usepackage{geometry}
\usepackage{mathtools}
\usepackage{color}
\usepackage[title]{appendix}
\usepackage[all]{xy}
\usepackage{tikz-cd}
\usepackage{blindtext}
\usepackage[T1]{fontenc}
\usepackage{amsthm}
\usepackage{amsfonts}
\usepackage{txfonts}
\usepackage{palatino,amssymb,epsfig}
\usepackage{latexsym,epsf,epic,amscd}
\usepackage{mathrsfs}
\usepackage{graphicx}
\usepackage{caption}
\usepackage{subcaption}
\usepackage{appendix}
\usepackage{amsmath}
\usepackage{bookmark}
\usepackage[nottoc]{tocbibind}
\hypersetup{
colorlinks=true,
linkcolor=black,
citecolor=blue,
anchorcolor=black,
urlcolor=black}
\allowdisplaybreaks[4]
\numberwithin{equation}{section}
\DeclareFontFamily{OMX}{yhex}{}
\DeclareFontShape{OMX}{yhex}{m}{n}{<->yhcmex10}{}
\DeclareSymbolFont{yhlargesymbols}{OMX}{yhex}{m}{n}
\DeclareMathAccent{\wideparen}{\mathord}{yhlargesymbols}{"F3}

\DeclareMathOperator{\area}{Area}

    \newcommand{\Addresses}{{
  \bigskip
  \footnotesize
  \noindent Te Ba, \href{batexu@hnu.edu.cn}{batexu@hnu.edu.cn}
  \newline\textit{School of Mathematics, Hunan University, Changsha 410082, P.R. China}\par\nopagebreak
  \medskip
  \noindent Guangming Hu, \href{18810692738@163.com}{18810692738@163.com}
 \newline\textit{ College of Science, Nanjing University of Posts and Telecommunications,
  Nanjing, 210003, P.R. China.}\par\nopagebreak
  \medskip
 \noindent Yu Sun, \href{yusun15185105160@163.com}{yusun15185105160@163.com}
  \newline\textit{School of Mathematics and Physics, Nanjing Institute of Technology, 211100, P.R. China.}\par\nopagebreak
}}
\title{Circle packings and hyperbolic surfaces of finite type}
\author{Te Ba,  Guangming Hu, Yu Sun}

\date{}
\newtheorem{theorem}{Theorem}[section]
\newtheorem{lemma}[theorem]{Lemma}

\theoremstyle{definition}

\newtheorem{remark}[theorem]{Remark}

\begin{document}
\maketitle

\begin{abstract}
This paper constructs hyperbolic polyhedral metrics via circle packings. We introduce the curvature of circles as a parameter to include all three types of constant curvature curves in the hyperbolic geometry. This provides a unified approach to producing polyhedral metrics for surfaces of broader topological types. The combinatorial total geodesic curvature serves as an effective tool for establishing the existence and uniqueness of the packing.

\end{abstract}

\section{Introduction}

\begin{figure}[htbp]
\centering
\captionsetup{width=0.90\linewidth}
\includegraphics[scale=0.22]{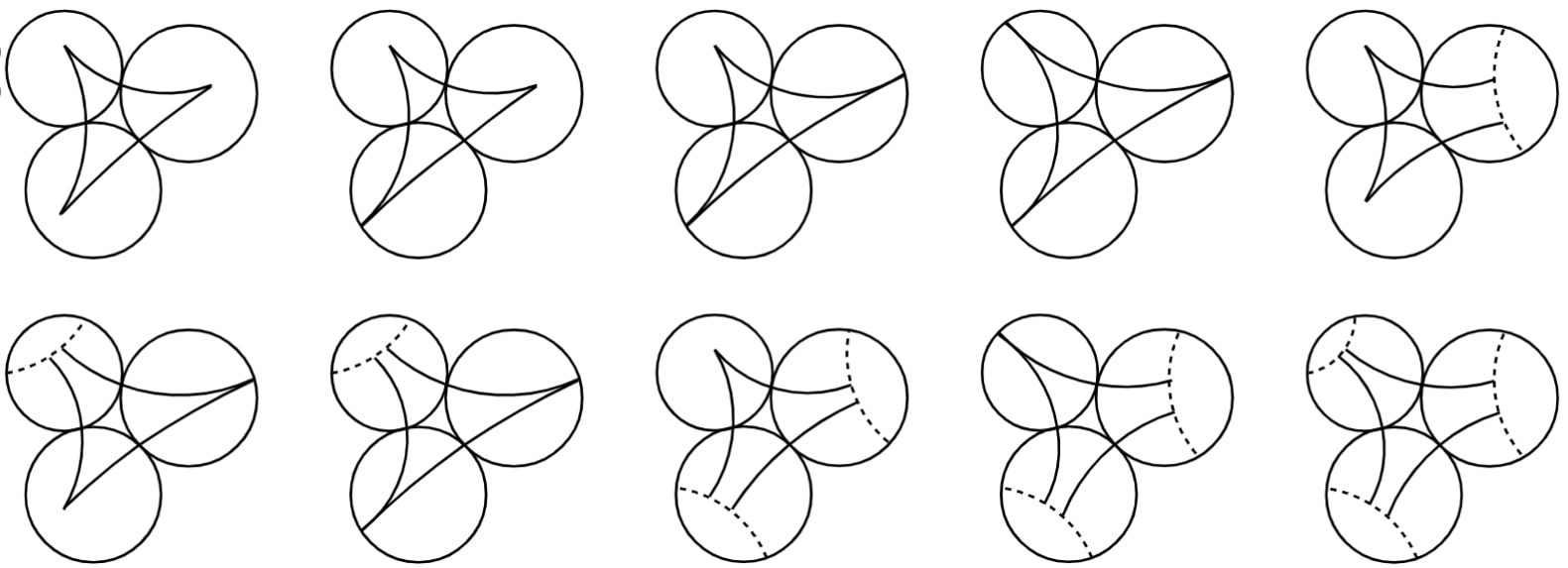}
\caption{Length assignment of decorated circle packings}
\label{F-1}
\end{figure}

Let $(\Sigma,\mathcal T)$ be a closed orientable $2$-manifold endowed with a triangulation $\mathcal{T}=\{V, E, F\}$.
Suppose $\kappa:V\to\mathbb{R}_+$ is a discrete function. We define
$V_P=\{i\in V\ |\ \kappa(i)=1\}$ and $V_\partial=\{i\in V\ |\ \kappa(i)<1\}.$
Let $\widetilde\Sigma \subset\Sigma$ be an embedded submanifold obtained by coning off all vertices in $V_P \subset V$ and removing a small open neighbourhood homeomorphic to a disk for each vertex in $V_\partial \subset V$.
One naturally obtains $(\widetilde\Sigma,\mathcal T)$, an orientable $2$-manifold $\widetilde\Sigma$ with cell decomposition induced by \(\mathcal{T}\). The face set is defined as $\widetilde F=\{uvw\cap\widetilde\Sigma\ |\ uvw\in F\}$,
and the edge set consists of $\widetilde E = \{uv \cap \widetilde\Sigma \mid uv \in E\}$ together with the boundary edges, i.e., those lying on $\partial \widetilde \Sigma$ and not contained in $\widetilde E$.
A \emph{decorated circle packing} on $(\Sigma,\mathcal T)$ determined by the discrete function $\kappa$ is a polyhedral metric on $(\widetilde\Sigma,\mathcal T)$ satisfying the following conditions:
\begin{itemize}
\item[$\textbf{(C1)}$] The length  $d(uv)$ of $uv\in\widetilde E$  is given by
\begin{equation}\label{E-1-1}
d(uv) =
\begin{cases}
\displaystyle \frac{1}{2} \ln\left( \frac{(\kappa(u) + 1)(\kappa(v) + 1)}{\left|(\kappa(u) - 1)(\kappa(v) - 1)\right|} \right), &   \kappa(u), \kappa(v) \neq 1, \\
+\infty, & \text{otherwise}.\\
\end{cases}
\end{equation}
\item[$\textbf{(C2)}$] Each interior angle at an endpoint of a boundary edge is equal to $\pi/2$.
\end{itemize}
Here, a hyperbolic polyhedral metric on  $(\widetilde\Sigma,\mathcal T)$ refers to a complete path metric that is piecewise hyperbolic on each face and the edge lengths coincide on adjacent faces.
The geometric meaning of (\ref{E-1-1}) is the distance between two (possibly ideal) centers or axes of a pair of mutually externally tangent curves in $\mathbb{H}^2$ with constant curvature $\kappa_u$, $\kappa_v$. Note that each constant curvature curve may be a circle, a horocycle, or a hypercycle.
The term \emph{decorated} is borrowed from the celebrated work of Penner~\cite{MR919235}, in which decorated horocycles are used to study the Teichm\"uller space of punctured surfaces.
Given the above settings, the polyhedral metric on each face belongs to one of the ten types listed in Figure~\ref{F-1}.

Decorated circle packings on surfaces related to several discrete geometric structures studied earlier.
We now present three cases that may be of particular interest to the reader.
 Set $V_{\mathrm int}=V\setminus (V_P\cup V_\partial)$.
A particularly classical case arises when $V_{\mathrm int}=V$,  so that $\widetilde{\Sigma}=\Sigma$ and the decorated circle packings reduce to the circle packings on closed triangulated surfaces.
This setting dates back to the work of Thurston~\cite{MR1435975} and have played significant roles in the studies of various areas including conformal mappings~\cite{MR906396,MR1207210}, hyperbolic polyhedra \cite{MR1283870,MR1370757,MR1943885,Zhou2023}, discrete conformal structures~\cite{MR3375525,MR3807319,MR3825607}, discrete analytic functions~\cite{MR2131318} combinatorics~\cite{MR1150601,MR3455159} and others. If $V_\partial=V$, then $\widetilde\Sigma$ is a compact surface with finite boundary components. A rich theory of discrete conformal structures on such surfaces has been studied in~\cite{MR2496046} and further extended in~\cite{MR2847231,MR4433084}. See the recent work~\cite{XuZheng2025} for a comprehensive survey of discrete conformal structures on such surfaces. Suppose $V_P=V$, making $\widetilde\Sigma$ a surface with $|V|$ punctures, a decorated circle packing on $(\Sigma,\mathcal T)$ constructs an ideal hyperbolic polyhedron, completely characterized in~\cite{MR1370757}. See Table~\ref{T-1} for a summary of face types in each case mentioned above. In all other cases, faces must admit more than one possible structure.

\begin{table}[ht]
\centering
\begin{tabular}{c|c}
\hline
 &  Face Configurations \\
\hline
$V_{\mathrm int}=V$ &  Triangles \\
\hline
$V_\partial=V$ &  Right-Angled Hexagons \\
\hline
$V_P=V$ &  Ideal Triangles \\
\hline
Otherwise &  Non-Uniform Face Structures \\
\hline
\end{tabular}
\caption{Three types of decorated circle packings on $(\Sigma,\mathcal T)$}
\label{T-1}
\end{table}
As stated in Theorem \ref{T-1-1}($i$), we show in this paper that any positive discrete function uniquely determines a decorated circle packing on surfaces.
This provides a new discrete framework that unifies the various geometric structures introduced above.
As a result of Theorem \ref{T-1-1}($i$) and the geometric meaning of $\textbf{(C1)}$, $\textbf{(C2)}$, a decorated circle packing on a surface gives rise to a collection of mutually externally tangent constant-curvature curves $\{\mathcal C_i\ |\ i\in V\}$, where each $\mathcal C_i$ has curvature $\kappa(i)$. We call each curve a decorated circle. Note that a decorated circle is piecewise smooth because singularities may exist at the center or along the axis of $\mathcal C_i$. The \emph{combinatorial total geodesic curvature} $L:V\to\mathbb{R}_+$ of a decorated circle packing is defined by
\[L(v)=T(C_v),\]
where $T(C_v)$ is the total geodesic curvature of the decorated circle $C_v$. See Figure \ref{F-2} for an illustration of three types of decorated circles on a vertex.

\begin{figure}[htbp]
  \centering
  \begin{minipage}[h]{0.3\textwidth}
    \centering
    \begin{tikzpicture}[scale=1.7, line join=round, line cap=round]
      \coordinate (A) at (0,0,0);
      \coordinate (B) at (2,0,0);
      \coordinate (C) at (1,1.4,0);
      \coordinate (D) at (1.7,-0.01,1);
      \tikzset{curve with mid/.style={decoration={markings,mark=at position 0.3 with \coordinate (#1);},postaction=decorate}}
      \path[draw=black, line width=0.7pt, bend right=10, curve with mid=A1] (A) to (C);
      \path[draw=black, line width=0.7pt, bend left=7, curve with mid=B1] (B) to (C);
      \path[draw=black, line width=0.7pt, bend left=10, curve with mid=D1] (D) to (C);
      \draw[line width=1.4pt, dashed, black] (A1) to[out=30, in=0, looseness=0.6] (B1);
      \draw[line width=1.4pt, black] (B1) to[out=10, in=-10, looseness=0.6] (D1);
      \draw[line width=1.4pt, black] (D1) to[out=-10, in=-50, looseness=0.6] (A1);
    \end{tikzpicture}
  \end{minipage}
  \begin{minipage}[h]{0.3\textwidth}
    \centering
    \begin{tikzpicture}[scale=1.7, line join=round, line cap=round]
      \coordinate (A) at (0,0,0);
      \coordinate (B) at (1.5,0,0);
      \coordinate (C) at (0.7,1.4,0);
      \coordinate (D) at (1.7,-0.01,1);
      \tikzset{curve with mid/.style={decoration={markings,mark=at position 0.3 with \coordinate (#1);},postaction=decorate}}
      \path[draw=black, line width=0.7pt, bend right=20, curve with mid=A1] (A) to (C);
      \path[draw=black, line width=0.7pt, bend left=20, curve with mid=B1] (B) to (C);
      \path[draw=black, line width=0.7pt, bend left=20, curve with mid=D1] (D) to (C);
      \draw[line width=1.4pt, dashed, black] (A1) to[out=30, in=0, looseness=0.6] (B1);
      \draw[line width=1.4pt, black] (B1) to[out=5, in=-10, looseness=0.6] (D1);
      \draw[line width=1.4pt, black] (D1) to[out=-20, in=-50, looseness=0.6] (A1);
    \end{tikzpicture}
  \end{minipage}
  \begin{minipage}[h]{0.33\textwidth}
    \centering
    \begin{tikzpicture}[scale=1.7, line join=round, line cap=round]
      \coordinate (A) at (0,0,0);
      \coordinate (B) at (1.7,0,0);
      \coordinate (D) at (1.5,-0.01,1);
      \coordinate (A2) at (0.5,1.3,0);
      \coordinate (B2) at (1.3,1.3,0);
      \coordinate (D2) at (1,1,0);
      \tikzset{curve with mid/.style={decoration={markings,mark=at position 0.2 with \coordinate (#1);},postaction=decorate}}
      \path[draw=black, line width=0.7pt, bend right=10, curve with mid=A1] (A) to (A2);
      \path[draw=black, line width=0.7pt, bend left=7, curve with mid=B1] (B) to (B2);
      \path[draw=black, line width=0.7pt, bend left=10, curve with mid=D1] (D) to (D2);
      \draw[line width=0.7pt, black] (A2) to[out=30, in=0, looseness=0.5] (B2);
      \draw[line width=0.7pt, black] (B2) to[out=10, in=-10, looseness=0.5] (D2);
      \draw[line width=0.7pt, black] (D2) to[out=-10, in=-50, looseness=0.5] (A2);
      \draw[line width=1.5pt, dashed, black] (A1) to[out=30, in=0, looseness=0.6] (B1);
      \draw[line width=1.5pt, black] (B1) to[out=10, in=-10, looseness=0.6] (D1);
      \draw[line width=1.5pt, black] (D1) to[out=-10, in=-50, looseness=0.6] (A1);
    \end{tikzpicture}
  \end{minipage}
  \caption{Three types of decorated circles highlighted with thick lines}
  \label{F-2}
\end{figure}

The notion of combinatorial total geodesic curvature was first introduced in the work of Nie~\cite{MR4683863} to generalize the existence and rigidity results for Delaunay circle patterns established by Bobenko and Springborn~\cite{MR2022715} in the spherical setting. Surprisingly, we find that combinatorial total geodesic curvature serves as an effective geometric data for studying the existence and uniqueness of decorated circle packings.

The main result of the paper is as follows.
\begin{theorem}\label{T-1-1}
Let $(\Sigma,\mathcal T)$ be a closed orientable triangulated 2-manifold. Then the following two statements hold:
 \begin{itemize}
 \item[($i$)] For any discrete function $\kappa:V\to\mathbb{R}_+$,  there exists a decorated circle packing on  $(\Sigma,\mathcal T)$ determined by $\kappa$, which is unique up to isometry.
     \item[($ii$)] A decorated circle packing on  $(\Sigma,\mathcal T)$ with combinatorial total geodesic curvature $L:V\to\mathbb{R}_+$ exists if and only if
          \[\sum_{i\in I} L(i)<\pi\vert F_I\vert\]
for every $I\subset V$, where $F_I$ denotes the set of faces having at least one vertex in $I$.  Moreover, such a decorated circle packing is unique up to isometry if it exists.
 \end{itemize}
\end{theorem}
\begin{remark}
Angle deficits $\Theta(v)$ at each vertex have been used as an important discrete curvature to study the rigidity of circle patterns in Euclidean and hyperbolic background geometries in previous works~\cite{MR1435975,MR3810260,MR4334399}. By comparison, the combinatorial total geodesic curvature of a Euclidean circle pattern at $v\in V$ can be calculated by
\[L(v)=l(v)\kappa(v)=\Theta(v)r(v)\frac{1}{r(v)}=\Theta(v),\]
where $l(v)$ and $\kappa(v)$ represent the length and the curvature of the circle at $v$. This coincides with the angle deficit.
However, the angle deficit fails to apply to decorated circle packings, as it is always zero at vertices in \(V_P\) and undefined at vertices in \(V_\partial\) due to the absence of interior angle structures.
\end{remark}
\textbf{Organization.} Section~2 presents a detailed analysis of three mutually tangent curves of constant curvature in $\mathbb{H}^2$ and provides a proof of Theorem~\ref{T-1-1}($i$). The existence argument of Theorem~\ref{T-1-1}($ii$) is established by characterizing the image of the curvature map, while the uniqueness part is proved using variational principles. These will be developed in Section~3. As an additional advantage, the variational principle offers an effective method for finding the desired decorated circle packings with prescribed combinatorial total geodesic curvature. This method will be demonstrated in Section~4 and formulated as Theorem~\ref{T-5-1}.

\section{Tangency configurations}
This section is dedicated to a detailed analysis on externally tangent configurations of three constant-curvature curves in $\mathbb{H}^2$, corresponding to the three-circle patterns associated with a single face in classical circle packings.
Throughout the section, two curves are said to be externally tangent if, in a neighborhood of the point of tangency, they lie in distinct regions divided by the common tangent geodesic.
The following lemma is due to the second author et al~\cite{MR4878836}. We present a proof for the sake of completeness.

\begin{lemma}\label{L-2-1}
For any $\kappa_1, \kappa_2, \kappa_3 > 0$, there exists a triple of smooth arcs $\mathcal C_i: [0,1]\to\mathbb{H}^2$ ($i=1,2,3$) such that each $\mathcal C_i$ has constant curvature $\kappa_i$ and each pair of arcs is externally tangent at their endpoint. Moreover, this triple is unique up to isometries of $\mathbb{H}^2$.
\end{lemma}

\begin{proof}
It can be verified that exactly three points of tangency are required to satisfy the condition of external tangency.
Let $\kappa=\max\{\kappa_1,\kappa_2,\kappa_3\}$. We divide the proof into three cases.
\begin{itemize}
\item[($a$)] Suppose $\kappa>1$. We will construct $\mathcal C_1, \mathcal C_2, \mathcal C_3$ in the Poincar\'e disk model $\mathbb{H}^2$. Without loss of generality, we assume that $\kappa_1>1$. Let us define
 \[\mathcal C_1=\big\{r_1e^{i\theta}\mid\theta\in[0,2\pi]\big\}\]
 and $A:=r_1$, a point in the Poincar\'e disk, where \( r_1 = \coth^{-1} \kappa_1 \) denotes the radius of \( \mathcal{C}_1 \). Then we can find two curves externally tangent to $\mathcal C_1$ at $A$ with curvature $\kappa_2$ and $\kappa_3$, denoted by $\mathcal C_2$ and $\gamma$, respectively. By external tangency, the Euclidean centers of $\mathcal C_2$ and $\gamma$ necessarily lie inside the unit circle. Let $\phi_{\alpha}:\mathbb{H}^2\to\mathbb{H}^2$ be the isometry given by $\phi_{\alpha}(z)=e^{i\alpha}z$. We can find exactly two $\beta_1,\beta_2\in(-\pi,\pi)$ with $\beta_2=-\beta_1$ such that $\phi_{\beta_i}(\gamma)$ is externally tangent to $\mathcal C_2$ on $\mathbb{R}^2$. The point of tangency lies at the midpoint of the Euclidean line segment joining the centers of $\mathcal C_2$ and $\phi_{\beta_i}(\gamma)$, and is therefore contained in $\mathbb{H}^2$. Note that $\phi_{\beta_i}(\gamma)$ is also externally tangent to $\mathcal C_1$ and has curvature $\kappa_3$. Thus we can denote $\mathcal C_3=\phi_{\beta_1}(\gamma)$.
 Then we restrict each \(\mathcal{C}_s\) to the subarc between its two points of tangency, which yields the required arcs.
The uniqueness follows from the fact that the reflection across the real axis brings $\phi_{\beta_1}(\gamma)$ to $\phi_{\beta_2}(\gamma)$.
\item[($b$)] Suppose $\kappa=1$. We will construct $\mathcal C_1, \mathcal C_2, \mathcal C_3$ in the upper half-plane model $\mathbb{H}^2$. Assume that $\kappa_1=1$. We define
    \[\mathcal C_1=\{x+i\ |\ x\in\mathbb{R}\}\]
    and set the point $B=i$. Then there exist two curves externally tangent to $\mathcal C_1$ at $B$ with constant curvatures $\kappa_2$ and $\kappa_3$, denoted by $\mathcal C_2$ and $\zeta$.
    The external tangency ensures that their Euclidean centers lie in $\mathbb{H}^2$.
    For $a\in\mathbb{R}$, we define an isometry $\psi_a:\mathbb{H}^2\to\mathbb{H}^2$ given by $\psi_a(z)=z+a$.
    We can find exactly two $b_1, b_2\in\mathbb{R}$ such that $\psi_{b_i}(\zeta)$ is externally tangent to $\mathcal C_1$ and $\mathcal C_2$. One can select $\mathcal C_3=\psi_{b_1}(\zeta)$ and make a restriction similar to that in case $(a)$ to derive the existence.  The reflection across the imaginary axis maps $\psi_{b_1}(\zeta)$ to $\psi_{b_2}(\zeta)$, which implies the uniqueness.

\item[($c$)] Suppose $\kappa<1$.  There exists a unique right-angled hyperbolic hexagon whose three non-pairwise adjacent edge lengths are
    \[l_w=\tanh^{-1}\kappa_u+\tanh^{-1}\kappa_v,\]
    where $\{u,v,w\}=\{1,2,3\}$.
    Then one can construct an equidistant curve at orthogonal distance $\tanh^{-1}\kappa_i$ from the side opposite to the side of length $l_s$ for $s=1,2,3$, denoted by $\mathcal C_s$, respectively. Note that the curvature of $\mathcal C_s$ is $\kappa_s$ and we say that two sides of a hexagon are opposite to each other if they are separated by exactly two other sides. Then $\mathcal C_u$ and $\mathcal C_v$ are externally tangent at some point on the side of length $l_w$.

\end{itemize}

\end{proof}
With Lemma~\ref{L-2-1} at hand, we are ready to prove Theorem~\ref{T-1-1}($i$).
Lemma~\ref{L-2-1} provides a unified structure of each face under a decorated circle packing, allowing us to proceed without resorting to a case-by-case analysis.
\begin{proof}[\textbf{Proof of Theorem \ref{T-1-1}($i$)}]
Let $\kappa:V\to\mathbb{R}_+$ be a discrete function.
According to the definition of decorated circle packing, one naturally obtains $V_{\mathrm int}$,  $V_P$, $V_\partial$, $\widetilde\Sigma\subset\Sigma$ and a polyhedral metric on $(\widetilde\Sigma,\mathcal T)$ from the data of $\kappa$.
This establishes the existence.

Let $f$ be a face in $F$. One can find a unique triple of vertices $u,v,w\in V$ such that $f$ corresponds to the triangle formed by $u,v,w$ in $(\Sigma,\mathcal T)$.
Then any decorated circle packing restricted to $f$ induces three mutually externally tangent curves with curvatures $\kappa_u$, $\kappa_v$ and $\kappa_w$ on $f$. Using these curvatures, one can compute each edge length and angle not assigned by $\textbf{(C1)}$, $\textbf{(C2)}$. The uniqueness then follows directly from Lemma \ref{L-2-1}.
\end{proof}
To proceed with the analysis of tangency configurations, we introduce the following construction.
A \textbf{right-angled bigon} in $\mathbb{H}^2$ is a convex topological disc bounded by two smooth arcs $\mathcal{C}_1, \mathcal{C}_2 : [0,1] \to \mathbb{H}^2$ of constant curvature that share the same endpoints and intersect orthogonally at both ends. Here, we refer to $\mathcal{C}_1$ and $\mathcal{C}_2$ as the sides of the right-angled bigon.
One can check that a right-angled bigon always exists whenever the prescribed curvatures $\kappa_1$ and $\kappa_2$ of the two sides satisfy $\max\{\kappa_1, \kappa_2\} \geq 1$.
We refer the reader to the work of Lutz~\cite{MR4523869} for more background on this topic from the perspective of M\"obius transformations.

Suppose $\mathcal{C}_1$ and $\mathcal{C}_2$ are the two sides of a right-angled bigon with curvatures $\kappa_1$ and $\kappa_2$, respectively, where $\kappa_2 \geq 1$.
 Let us define the total geodesic curvature of $\mathcal{C}_1$ and $\mathcal{C}_2$ by
\begin{equation}\label{E-gh}g(\kappa_1,\kappa_2)=\int_{\mathcal C_1}\kappa_1ds,\quad h(\kappa_1,\kappa_2)=\int_{\mathcal C_2}\kappa_2ds.\end{equation}
We are fortunate to find that both \( g \) and \( h \) admit explicit expressions. For brevity, we denote $b_{\kappa_s}=\frac{1}{\sqrt{|\kappa_s^2-1|}}$. It turns out that
\begin{equation}\label{E-2-2}
g(\kappa_1,\kappa_2)
=
\begin{cases}
2b_{\kappa_1}\kappa_1\cot^{-1}(b_{\kappa_1}\kappa_2),
& \kappa_1>1,\\
\frac{2\kappa_1}{\kappa_2},
& \kappa_1=1,\\
2b_{\kappa_1}\kappa_1\coth^{-1}(b_{\kappa_1}\kappa_2),
& \kappa_1<1.
\end{cases}
\end{equation}
and
\begin{equation}\label{E-2-3}
h(\kappa_1,\kappa_2)
=2b_{\kappa_2}\kappa_2\cot^{-1}(b_{\kappa_2}\kappa_1).
\end{equation}
We now present a more refined result regarding these two functions. For the hyperbolic trigonometric identities used in the proof of Lemma \ref{L-3-3}, we refer the reader to the monograph of Buser~\cite[Chapter 2]{MR1183224} and the work of Guo and Luo~\cite[Appendix A]{MR2496046}.
A technical lemma used in the proof of Lemma~\ref{L-3-3}~($i$) is presented in the appendix, and the proof of Lemma~\ref{L-3-3}~($ii$) is also included in the appendix, as it involves a direct calculation.
\begin{figure}[htbp]
\centering
\includegraphics[scale=0.35]{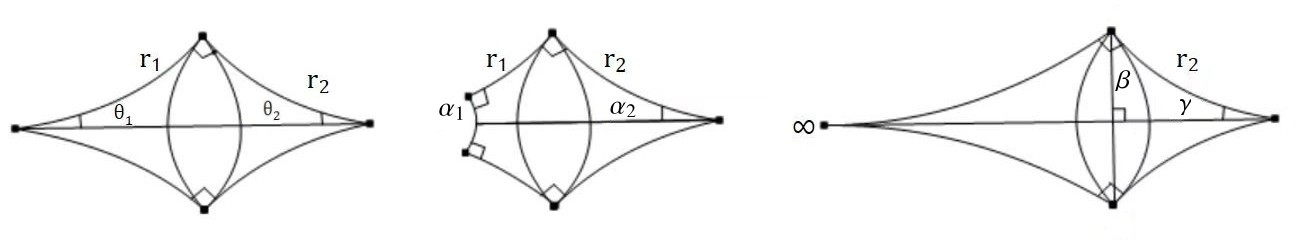}
\caption{Geometric relations on right-angle bigons}
\label{F-3}
\end{figure}

\begin{lemma}\label{L-3-3}
 Let $\mathcal C_1$ and $\mathcal C_2$ be two sides of a right-angled bigon in~$\mathbb H^2$ with curvatures $\kappa_1,$ and $\kappa_2$ respectively, and assume that $\kappa_2\geq 1$. Let $g$, $h$ be the total geodesic curvatures of $\mathcal C_1$ and $\mathcal C_2$ as defined in (\ref{E-gh}).
 Then the following statements hold:
\begin{itemize}
\item[($i$)] $g$, $h$ are given by (\ref{E-2-2}) and (\ref{E-2-3}).
\item[($ii$)]  $g$, $h$ are $C^1$-smooth functions.
\end{itemize}
\end{lemma}
\begin{proof}
In order to express $g$ and $h$, we proceed by considering the following three cases.
\begin{itemize}
\item[($a$)] Suppose $\kappa_1>1$. The centers of circles $\mathcal C_1$ and $\mathcal C_2$, together with their intersection point, form a right triangle, as illustrated in the left part of Figure~\ref{F-3}. We denote by $\theta_ s$  the interior angle of the right triangle at the center of $\mathcal C_s$ and by $r_s$ the radius of $\mathcal C_s$. Then
    \[g(\kappa_1,\kappa_2)=2\theta_1\cosh r_1=2\cot^{-1}(\coth r_2\sinh r_1)\cosh r_1,\]
which is equal to $h(\kappa_2,\kappa_1)$.
\item[($b$)] Suppose $\kappa_1<1$. Note that the axis of \(\mathcal C_1\) and the center of \(\mathcal C_2\), their point of intersection, and the perpendicular from the axis of \(\mathcal C_1\) and the center of \(\mathcal C_2\) together form a hyperbolic quadrilateral with three right angles, as shown in the middle part of Figure~\ref{F-3}. Let $\alpha_1$ denote the side length of the quadrilateral on the axis of $\mathcal C_1$ and $\alpha_2$ by the interior angle at the center of $\mathcal C_2$. By hyperbolic trigonometric identities, we obtain
    \[g(\kappa_1,\kappa_2)=2\alpha_1\sinh r_1=2\coth^{-1}(\coth r_2\cosh r_1)\sinh r_1,\]
    \[h(\kappa_1,\kappa_2)=2\alpha_2\cosh r_2=2\cot^{-1}(\sinh r_2\tanh r_1)\cosh r_2.\]
\item[($c$)] Suppose $\kappa_1=1$. The right triangle analogous to the first case becomes a right triangle with an ideal vertex. We obtain
    \begin{equation}\label{E-2-4}\sin\gamma_2\cosh r_2=1,\end{equation}
where $\gamma$ is the interior angle of the ideal right triangle at the center of $\mathcal C_2$.
    Let $\beta$ denote the angle between the geodesic connecting the point of intersection of  $\mathcal C_1$ and $\mathcal C_2$ and the line joining the point of intersection to the center of  $\mathcal C_2$, as illustrated in the right part of Figure~\ref{F-3}. Then
\[\tan\beta=\frac{\cot\gamma_2}{\cosh r_2}=\frac{\sinh r_2}{\cosh r_2}=\frac{1}{\kappa_2}.\]
By Lemma \ref{L-6-1}, we derive that
    \[g(\kappa_1,\kappa_2)=2 \kappa_1\tan\beta=\frac{2\kappa_1}{\kappa_2}.\]
    Also, we derive
    \[h(\kappa_1,\kappa_2)=2\gamma\cosh r_2=\cot^{-1}(\sinh r_2)\cosh r_2\]
from (\ref{E-2-4}). Then the explicit expressions of $g$ and $h$ follow directly from the following identities
\[\begin{aligned}
\kappa_s=\begin{cases}
\coth r_s, & \kappa_s>1,\\
\tanh r_s, & 0<\kappa_s<1.\\
\end{cases}
\end{aligned}\]
\end{itemize}

\end{proof}

The following lemma is from the second author et al~\cite{MR4878836}. With the explicit expressions of $g$ and $h$ at hand, we present an alternative proof, distinct from that in~\cite{MR4878836}.

\begin{lemma}\label{L-3-4}
Let $\mathcal{C}_u, \mathcal{C}_v, \mathcal{C}_w: [0,1] \to \mathbb{H}^2$ be three smooth arcs of constant curvatures $\kappa_u, \kappa_v, \kappa_w>0$, respectively. If each pair of them is externally tangent at a common endpoint, then there exists a unique constant curvature curve $\mathcal C$ passing through the three points of tangency and orthogonal to each $\mathcal C_i$. Moreover, the curvature $\kappa$ of $\mathcal C$ is given by
\begin{equation}\label{E-2-7}\kappa^2=\kappa_u\kappa_v+\kappa_u\kappa_w+\kappa_v\kappa_w+1.\end{equation}
\end{lemma}

\begin{proof}
Let \(P, Q, R\) be three points of tangency. By applying a suitable M\"obius transformation, we may assume without loss of generality that there exists a unique Euclidean circle \(\mathcal{C}\) passing through  \(P, Q, R\). Since a Euclidean circle represents a constant curvature curve in \(\mathbb{H}^2\), this proves the existence.

Let us compute the curvature of $\mathcal C$. We claim that $\mathcal C$ is orthogonal to each $C_i$. It follows from the fact that the Euclidean center of $\mathcal C$ is equidistant (in the Euclidean sense) from $P$, $Q$, $R$, which implies the center must lie on the three Euclidean tangent lines. Consequently, \( \mathcal{C} \) is orthogonal to each of the three arcs.
Recalling (\ref{E-2-3}), we can use
\begin{equation}\label{E-3-7}\theta(\kappa_s,\kappa)=\frac{h(\kappa_s,\kappa)}{\cosh r}=2\cot^{-1}(b_{\kappa_s}\kappa)\end{equation}
to represent the angle formed by the center of $\mathcal C$  and its two intersection points with $\mathcal C_1$ if $\kappa>1$, where $b_{\kappa_s}=\frac{1}{\sqrt{|\kappa_s^2-1|}}$ and $r=\coth^{-1}k$ is the radius of $\mathcal C$. This gives rise to a function $\Theta:(1,+\infty)\to\mathbb{R}$ defined by
\[\Theta(\kappa)=\theta(\kappa_u,\kappa)+\theta(\kappa_v,\kappa)+\theta(\kappa_w,\kappa)-2\pi.\]
Note that the geometric interpretation of $\Theta(\kappa)$ is the angle at the center of a circle of curvature $\kappa$, enclosed by three mutually tangent arcs of curvatures $\kappa_u$, $\kappa_v$, and $\kappa_w$, each intersecting the circle orthogonally.
Observe that $\theta(\kappa_s,\kappa)\to0$ as $\kappa\to1^+$ and $\theta(\kappa_s,\kappa)\to\pi$ as $\kappa\to+\infty$. It follows that $\Theta(\kappa)\to0$ as $\kappa\to1^+$ and $\Theta(\kappa)\to3\pi$ as $\kappa\to+\infty$. Furthermore, one can compute that
\begin{equation}\label{E-2-8}\frac{d\theta(\kappa_s,k)}{d\kappa}=\frac{2\kappa_s\kappa b_{\kappa}}{\kappa_s^2+\kappa^2-1}>0,\end{equation}
which implies $\Theta$ is strictly increasing.
By the Intermediate Value Theorem, one can find a unique $\kappa_0\in(1,+\infty)$ such that $\Theta(\kappa_0)=2\pi$.

To arrive at an explicit expression,  we make use of the equality
$$
\theta(\kappa_u,\kappa)+\theta(\kappa_v,\kappa)+\theta(\kappa_w,\kappa)=2\pi
$$
to deduce that
\[
\begin{aligned}
\cot\frac{\theta_u}{2}&=\cot(\pi-\frac{\theta_v}{2}-\frac{\theta_w}{2})\\
& =\frac{-\cos \frac{\theta_v}{2} \cos \frac{\theta_w}{2}+\sin \frac{\theta_v}{2} \sin \frac{\theta_w}{2}}{\sin \frac{\theta_v}{2} \cos \frac{\theta_w}{2}+\cos \frac{\theta_v}{2} \sin \frac{\theta_w}{2}}\\
&=\frac{-\cot \frac{\theta_v}{2} \cot \frac{\theta_w}{2}+1}{\cot \frac{\theta_v}{2}+\cot \frac{\theta_w}{2}},\\
\end{aligned}
\]
where we write $\theta_s=\theta(\kappa_s,\kappa)$ for short.
Substituting (\ref{E-3-7}) into the above equation gives \eqref{E-2-7}.

\end{proof}

Equipped with these two lemmas, we proceed to a more detailed analysis of the tangency configurations.
Let $\mathcal C_u$, $\mathcal C_v$ and $\mathcal C_w$ be three smooth arcs of curvature  $\kappa_u$, $\kappa_v$ and $\kappa_w$, respectively, and each pair of arcs is externally tangent at their endpoint.
The total geodesic curvature function of $\mathcal C_i$ is given by
\begin{equation}\label{E-3-1}f_i(\kappa_u,\kappa_v,\kappa_w)=\int_{\mathcal C_i}\kappa_ids\end{equation}
for $i\in\{u,v,w\}$. Under substitution $k_i=\ln\kappa_i$,  it gives rise to the following differential $1$-form
\[
\omega=f_udk_u+f_vdk_v+f_wdk_w.
\]
Furthermore, we will establish the following refined result for tangency configurations.
\begin{lemma}\label{L-3-1}
Let $f_i$ be defined as in (\ref{E-3-1}). Then the following statements hold:
\begin{itemize}
\item[($i$)]  Each $f_i$ is $C^1$-smooth;
\item[($ii$)] The $1$-form $\omega$ is closed.
\end{itemize}
\end{lemma}

\begin{proof}
First, we prove ($i$). By Lemma \ref{L-3-4}, there exists a unique constant curvature curve $\mathcal C$ passing through all three points of tangency and orthogonal to $\mathcal C_u$, $\mathcal C_v$ and $\mathcal C_w$. Then $\mathcal C$ forms a right-angled bigon with each $\mathcal C_i$. We use $\kappa_c$ to denote the curvature of $\mathcal C$.
For any $\kappa_u,\kappa_v,\kappa_w>0,$ we know
\[f_i(\kappa_u,\kappa_v,\kappa_w)=g(\kappa_i,\kappa_c)\]
for $i\in\{u,v,w\}$ since $\kappa_c>1$.
Then the $C^1$-smoothness of $f_i$ follows from $C^1$-smoothness of $g$.

Then we turn to the proof of $(ii)$. Since $g$ is $C^1$, a straightforward calculation based on (\ref{E-2-2}) yields that
\begin{equation}\label{E-2-9}\
g_y(\kappa_1,\kappa_2)=\frac{2\kappa_1}{1-\kappa_1^2-\kappa_2^2}.
\end{equation}
Combining (\ref{E-2-7}) and (\ref{E-2-9}), we derive that
\begin{equation}\label{E-2-10}\begin{aligned}
\frac{\partial}{\partial k_v}f_u(\kappa_u,\kappa_v,\kappa_w)&=g_y(\kappa_u,\kappa_c)\frac{\partial \kappa_c}{\partial \kappa_v}\frac{d\kappa_v}{dk_v}&=\frac{2\kappa_u\kappa_v}{1-\kappa_u^2-\kappa_c^2}\frac{\kappa_u+\kappa_w}{2\kappa_c},\\
\frac{\partial}{\partial k_u}f_v(\kappa_u,\kappa_v,\kappa_w)&=g_y(\kappa_v,\kappa_c)\frac{\partial \kappa_c}{\partial \kappa_u}\frac{d\kappa_u}{dk_u}&=\frac{2\kappa_u\kappa_v}{1-\kappa_v^2-\kappa_c^2}\frac{\kappa_v+\kappa_w}{2\kappa_c}.
\end{aligned}\end{equation}
By Lemma \ref{L-3-4}, we deduce that 
\[\begin{aligned}
(\kappa_u+\kappa_w)(1-\kappa_v^2-\kappa_c^2)-(\kappa_v+\kappa_w)(1-\kappa_u^2-\kappa_c^2)
&=(\kappa_u-\kappa_v)(\kappa_u\kappa_v+1)-(\kappa_u-\kappa_v)\kappa_c^2+(\kappa^2_u-\kappa^2_v)\kappa_w\\
&=(\kappa_u-\kappa_v)(\kappa_u\kappa_v+\kappa_u\kappa_w+\kappa_v\kappa_w-\kappa_c^2+1)\\
&=0.
\end{aligned}\]
Then we immediately obtain \[\frac{\partial}{\partial k_v}f_u(\kappa_u,\kappa_v,\kappa_w)=\frac{\partial}{\partial k_u}f_v(\kappa_u,\kappa_v,\kappa_w).\]
\end{proof}

The following lemma describes two key monotonicity properties of the total geodesic curvature functions associated with a tangency configuration.
\begin{lemma}\label{L-3-5}
Let $\mathcal C_u,\mathcal C_v,\mathcal C_w: [0,1]\to\mathbb{H}^2$ be three smooth arcs of constant curvatures $\kappa_u, \kappa_v, \kappa_w>0$, respectively, such that each pair is externally tangent at their common endpoint.
Let $f_i$ be defined as in (\ref{E-3-1}). Then
\begin{subequations}
\begin{align}
\label{E-3-9a}\frac{\partial f_u}{\partial \kappa_v}&<0,\\
\label{E-3-9b}\frac{\partial(f_u+f_v+f_w)}{\partial \kappa_u}&>0.
 \end{align}\end{subequations}
\end{lemma}
\begin{proof}
We readily obtain (\ref{E-3-9a}) from (\ref{E-2-10}). Let us prove (\ref{E-3-9b}).
We use $\kappa_c$ to denote the curvature passing through three points of tangency.
To simplify the computations, we derive the substitution
\begin{equation}\label{E-3-10}
\frac{\partial \kappa_c}{\partial \kappa_u}=\frac{-\theta_x(k_u,k_c)}{\theta_y(k_u,k_c)+\theta_y(k_v,k_c)+\theta_y(k_w,k_c)},
\end{equation}
which arises from differentiating
 \[\theta(k_u,k_c)+\theta(k_v,k_c)+\theta(k_w,k_c)=2\pi\]
on both sides with respect to $\kappa_u$. Meanwhile, we compute from (\ref{E-3-7}) that
\begin{equation}\label{E-3-3}\theta_x(\kappa_u,\kappa_c)=\frac{2}{b_{\kappa_c}(1-\kappa_u^2-\kappa_c^2)}.\end{equation}
Comparing (\ref{E-2-8}) and (\ref{E-2-9}), we observe that
\begin{equation}\label{E-3-4}\theta_y(\kappa_i,\kappa_c)=-g_y(\kappa_i,\kappa_c)\kappa_cb_{\kappa_c}\end{equation}
for $i\in\{u,v,w\}$. According to (\ref{E-3-10}), (\ref{E-3-3}) and (\ref{E-3-4}), we deduce that
 \[\begin{aligned}
 \frac{\partial}{\partial \kappa_u}(f_u+f_v+f_w)&=\frac{\partial}{\partial k_u}(g(\kappa_u,\kappa_c)+g(\kappa_v,\kappa_c)+g(\kappa_w,\kappa_c))\frac{\partial \kappa_c}{\partial \kappa_u}\\
 &=g_x(\kappa_u,\kappa_c)+\big(g_y(\kappa_u,\kappa_c)+g_y(\kappa_v,\kappa_c)+g_y(\kappa_w,\kappa_c)\big)\frac{\partial \kappa_c}{\partial \kappa_u}\\
 &=g_x(\kappa_u,\kappa_c)+\frac{\theta_x(\kappa_u,\kappa_c)}{\kappa_cb_{\kappa_c}}\\
&=g_x(\kappa_u,\kappa_c)+\frac{2}{\kappa_c(1 - \kappa_u^2 - \kappa_c^2)b^2_{\kappa_c}}.
 \end{aligned}\]
We now apply the result of  Lemma \ref{L-3-3}($ii$), (\ref{E-6-1}) and (\ref{E-6-3}) to obtain

  \[
\frac{\partial}{\partial k_u}(f_u+f_v+f_w)
=
\begin{cases}
\frac{2}{(k_u^2-1)k_c}-2b^3_{k_u}\cot^{-1}b_{k_u}k_c,
& k_u>1,\\
\frac{2}{3k_c^3}
& k_u=1,\\
\frac{2}{(k_u^2-1)k_c}+2b^3_{k_u}\cot^{-1}b_{k_u}k_c,
& k_u<1.
\end{cases}
\]
Observe that $\cot^{-1}b_{k_u}k_c<\frac{1}{b_{k_u}k_c}$ when $k_u\neq 1$. It follows that
\[\frac{\partial}{\partial k_u}(f_u+f_v+f_w)>0.\]
\end{proof}

The next lemma characterizes the limit behavior of $f_u$, $f_v$ and $f_w$.

\begin{lemma}
Let $\mathcal C_u,\mathcal C_v,\mathcal C_w: [0,1]\to\mathbb{H}^2$ be three smooth arcs of constant curvatures $\kappa_u$, $\kappa_v$ and $\kappa_w$, respectively, such that each pair is externally tangent at their common endpoint. Let $0\leq a,b,c<+\infty$. Then the following statements hold:
\begin{subequations}
\begin{align}
    \label{E-3-3a}\lim_{\kappa_u\to 0}f_u&=0,\\
\label{E-3-3b}\lim_{(\kappa_u,\kappa_v,\kappa_w)\to(+\infty,a,b)}f_u&=\pi,\\
\label{E-3-3c}\lim_{(\kappa_u,\kappa_v,\kappa_w)\to(+\infty,+\infty,c)}f_u+f_v&=\pi,\\
\label{E-3-3d}\lim_{(\kappa_u,\kappa_v,\kappa_w)\to(+\infty,+\infty,+\infty)}f_u+f_v+f_w&=\pi.
\end{align}\end{subequations}

\end{lemma}

\begin{proof}
We divide the proof of (\ref{E-3-3a}) into two cases.
\begin{itemize}
\item[($i$)] Suppose the limits of $\kappa_v$ and $\kappa_w$ equals to $0$.
Then $\area(\Lambda)\to\pi$ as $\Lambda$ approaches to an ideal triangle. By Lemma \ref{L-3-7}, we obtain \[f_u+f_v+f_w\to0,\]
which implies $f_u\to0$.
\item[($ii$)] Suppose the limit of at least one of $\kappa_v,\kappa_w$ is greater than $0$. Let $\{(\kappa_u^m,\kappa_v^m,\kappa_w^m)\}_{m\in\mathbb{N}}$ be a sequence with limit $(0,\lambda,\mu)$. Let us consider the sequence $\{(\kappa_u^m,\frac{\kappa_v^m}{m},\frac{\kappa_w^m}{m})\}_{m\in\mathbb{N}}$. By (\ref{E-3-9a}), we deduce that
\[
f_u(\kappa_u^m,\kappa_v^m,\kappa_w^m)\leq f_u(\kappa_u^m,\frac{\kappa_v^m}{m},\frac{\kappa_w^m}{m}).
\]
Then the proof reduces to the first case.
\end{itemize}
We turn to prove (\ref{E-3-3b}).
Let us denote the length of $\mathcal C_s$ by $l_s$. Note that
\[f_u=\pi-\area(\Lambda)-f_v-f_w.\]
We observe that $f_v,f_w\to 0$ as $l_y,l_w\to0$. From Lemma \ref{L-3-4}, we know $\kappa_c\to+\infty$ as $\kappa_u\to+\infty$, where $k_c$ represents the curvature of the circle passing three points of tangency. This leads to $\area(\Lambda)\to0$, which yields $f_u\to\pi$. The proof of (\ref{E-3-3c}), (\ref{E-3-3d}) are analogous to the proof of (\ref{E-3-3b}).

\end{proof}
The following lemma is a direct result of the Gauss-Bonnet Theorem.
\begin{lemma}\label{L-3-7}
Let $\mathcal C_u,\mathcal C_v,\mathcal C_w: [0,1]\to\mathbb{H}^2$ be three smooth arcs of constant curvatures $\kappa_u$, $\kappa_v$ and $\kappa_w$, respectively, such that each pair is externally tangent at their common endpoint. Let $\Lambda$ be the region enclosed by $\mathcal C_u,\mathcal C_v$ and $\mathcal C_w$ and let $f_i$ be defined as in \eqref{E-3-1}. Then
\[\area(\Lambda)=\pi-f_i-f_j-f_k.\]
\end{lemma}

\section{Image of Curvature Maps}
In this section, we shall prove Theorem \ref{T-1-1}($ii$).
The existence argument relies on a continuity method, originating from the foundational work of Thurston~\cite{MR1435975},  and was further elaborated by Marden and Rodin~\cite{MR1071766} and others.
The proof of uniqueness employs the variational principle introduced by Colin de Verdi\`ere~\cite{MR1106755}, which has since become an important tool for the study of various discrete geometric structures; see~\cite{MR2233848,MR2788656,MR2862158} and others.

Based on Theorem \ref{T-1-1}($i$) established in Section~2, any positive discrete function on $V$ naturally corresponds to a decorated circle packing on $(\Sigma,\mathcal T)$.
Let $\kappa:V\to\mathbb{R}$ be such a decorated circle packing on $(\Sigma,\mathcal T)$.
Recall that \( (\Sigma, \mathcal{T}) \) admits a collection of mutually externally tangent curves \( \{\mathcal{C}_v \mid v \in V\} \), each having curvature \( \kappa(v) \), and the combinatorial total geodesic curvature of the decorated circle packing at a vertex \( v \) is defined as the total geodesic curvature of \( \mathcal{C}_v \).
By (\ref{E-3-1}), the combinatorial total geodesic curvature at $u\in V$ can be computed by
\[
L_u=\sum_{uvw\in F} f_u(\kappa_u,\kappa_v,\kappa_w),
\]
 where we write $\kappa(u)=\kappa_u$ and denote a face by a triple $uvw$ for brevity.
By assigning labels $1,\cdots,n$ to all vertices in $V$, we define the differential $1$-form
\[
\eta=\sum_{i=1}^{|V|}L_idk_i
\]
on $\mathbb{R}^{|V|}$ under the coordinate change \(k_s = \ln \kappa_s\). Based on Lemma \ref{L-3-1}, one can verify that $\eta$ is closed because
\[\frac{\partial L_u}{\partial k_v}=\frac{\partial f_u(\kappa_u,\kappa_v,\kappa_w)}{\partial k_v}=\frac{\partial f_v(\kappa_u,\kappa_v,\kappa_w)}{\partial k_u}=\frac{\partial L_v}{\partial k_u}\]
if $uv\in E$ and $\frac{\partial L_u}{\partial k_v}=0$ otherwise.
Then we can define the function $\Phi:\mathbb{R}^{|V|}\to\mathbb{R}$ by
\begin{equation}\label{E-4-1}
\Phi(k)=\int_{k_0}^k\eta.
\end{equation}
\begin{lemma}\label{L-4-1}
The function $\Phi$ is strictly convex.
\end{lemma}
\begin{proof}
It suffices to show that the Hessian of $\Phi$ is positive definite. By Lemma \ref{L-3-5}, we know
\[\frac{\partial^2\Phi}{\partial k_u\partial k_v}=\frac{\partial L_u}{\partial k_v}=\sum_{uvw\in F}\frac{\partial}{\partial \kappa_v} f_u(k_u,k_v,k_w)\frac{d\kappa_v}{dk_v} <0\]
if $uv\in E$ and equals zero, otherwise.
Additionally, we have

\[\frac{\partial^2\Phi}{\partial k_v^2}=\frac{\partial L_v}{\partial k_v}=\sum_{uvw\in F}\frac{\partial}{\partial \kappa_v} f_v(k_u,k_v,k_w)\frac{d\kappa_v}{dk_v}>0,\]
which can be derived directly by combining (\ref{E-3-9a}) and (\ref{E-3-9b}). It follows that

\[\begin{aligned}
\left|\frac{\partial^2\Phi}{\partial k_v^2}\right|-\sum_{u\neq v}\left|\frac{\partial^2\Phi}{\partial k_u\partial k_v}\right|&=\frac{\partial }{\partial k_v}\bigg(L_v+\sum_{uv\in E}L_u\bigg)\\
&=\frac{\partial }{\partial \kappa_v}\bigg(\sum_{uvw\in F} f_u(k_u,k_v,k_w)+ f_v(k_u,k_v,k_w)+f_w(k_u,k_v,k_w)\bigg)\frac{d\kappa_v}{dk_v}\\
&>0,
\end{aligned}\]
which follows from (\ref{E-3-9b}). This indicates that the Hessian of $\Phi$ is diagonally dominant, which implies that it is positive definite.
\end{proof}
We consider the gradient map of  $\Phi$ given by
\begin{equation}
\begin{aligned}\label{E-4-2}
\nabla\Phi\;:\;\mathbb{R}^{|V|}\;&\longrightarrow\;\mathbb{R}^{|V|}\\
\tau\;\;\;\,&\longmapsto\;\;\nabla\Phi(\tau).
\end{aligned}
\end{equation}
By the definition of energy function $\Phi$, we observe that the gradient map of  $\Phi$ sends a decorated circle packing on $(\Sigma,\mathcal T)$ to its associated combinatorial total geodesic curvature.
Furthermore, we derive the following result.

\begin{lemma}\label{L-4-2}
Let $\Phi$ be defined as (\ref{E-4-1}). Then the gradient map $\nabla\Phi$ is injective and image of $\nabla\Phi$ consists of vectors $(L_1,\cdots,L_{|V|})\in\mathbb{R}^{|V|}_{+}$ satisfying
\[\sum_{i\in I}L_i<\pi\vert F_I\vert\ \text{for any subset}\ I\subset V\]
for each $I\subset V$,  where $F_I$ is the set of faces having at least one vertex in $I$ for the subset $I\subset V$.
\end{lemma}
\begin{proof}
The injectivity follows from Lemma \ref{L-4-1} and a classical result of analysis, which is stated as Lemma \ref{L-4-3}. Suppose $(L_1,\cdots,L_{|V|})\in\mathrm{Im}(\nabla\Phi)$.
For a subset $I\subset V$, we observe that
\[\sum_{u\in I}L_u=\sum_{u\in I}\sum_{uvw\in F}f_u(k_u,k_v,k_w)<\sum_{u\in I}\sum_{uvw\in F}\pi=|F_I|\pi,\]
where the inequality follows from Lemma \ref{E-3-7}.
By Brouwer's Theorem on the Invariance of Domain, it suffices to characterize the boundary of $\nabla\Phi(\mathbb{R}^{|V|})$ by analyzing the limit behavior of $\nabla\Phi(k^{(m)})$ for $\{k^{(m)}\}_{m\in\mathbb{N}}\subset\mathbb{R}^{|V|}$ as $\Vert k^{(m)}\Vert\to+\infty$.
For brevity, we denote 
\[k^{(m)}=(k_1^{(m)},\cdots,k_{|V|}^{(m)}).\]
Let $U_-,U_+\subset V$ be two subset such that $k_v^{(m)}\to-\infty$ for $v\in U_-$ and $k_v^{(m)}\to+\infty$ for $v\in U_+$. For each $v\in U_-$, we know
\[L_v^{(m)}=\sum_{uvw\in F}f_v(k_u^{(m)},k_v^{(m)},k_w^{(m)})\to 0\]
 as each term approaches zero, which follows from (\ref{E-3-3a}). We denote $F_i\subset F$ by the set of face having $i$ vertices in $U_+$ and we write $f_u(k_u^{(m)},k_v^{(m)},k_w^{(m)})=f^{(m)}_u$ for short. It follows that
\[\begin{aligned}
\sum_{v\in U_+}L_v^{(m)}&=\sum_{v\in U_+}\sum_{uvw\in F}f_v(k_u^{(m)},k_v^{(m)},k_w^{(m)}),\\
&=\sum_{f\in F_1}f^{(m)}_u+\sum_{f\in F_2}f^{(m)}_u+f^{(m)}_v+\sum_{f\in F_3}f^{(m)}_u+f^{(m)}_v+f^{(m)}_w\\
&\to |F_{U_2}|\pi
\end{aligned}\]
as the limit of each summation approach to $\pi$, which follows from (\ref{E-3-3b})-(\ref{E-3-3d}). Then the limit of $\nabla\Phi(k^{(m)})$ lies on the intersection of the planes which given by $L_i=0$ for each $i\in U_-$ and
\[\sum_{i\in U_+}L_i=\pi\vert F_{U_+
}\vert.\]
This completes the proof, as the sequence was chosen arbitrarily.
\end{proof}
\begin{lemma}\label{L-4-3}
Suppose $\Omega$ is convex and the $C^2$-smooth function $h:\Omega\to\mathbb{R}$ is strictly convex. Then the gradient map $\nabla h:\Omega\to\mathbb{R}^n$ is injective.
\end{lemma}

\begin{proof}[\textbf{Proof of Theorem \ref{T-1-1}($ii$)}]
As established in Lemma~\ref{L-3-1}, the total geodesic curvature function $f_i$ is $C^1$-smooth, which implies that the energy function $\Phi$ is $C^2$-smooth. Then the uniqueness argument directly follows from Lemma \ref{L-4-3}, while existence is an immediate consequence of Lemma~\ref{L-4-2}.

\end{proof}

\section{Combinatorial Curvature flows}

In this section, we provide an algorithm for finding the decorated circle packing with prescribed combinatorial total geodesic curvature.
For this purpose, we follow the work of Chow and Luo~\cite{MR2015261} on combinatorial Ricci flows, which serves as a guiding framework for several combinatorial geometric flows~\cite{MR2100762,MR3729504,MR4024610,MR3914484,MR4466650} and others.

Recall that the gradient of $\Phi$ gives an injective map from a decorated circle packing to its combinatorial total geodesic curvatures.
We define a function $\mathcal E:\mathbb{R}^{|V|}\to\mathbb{R}$ given by
\[\mathcal E=\Phi-\sum_{v\in V}\hat{L_v}k_v,\]
where each $\hat{L_v}$ is a given real constant. We then consider negative gradient flow of $\mathcal E$ governed by
\begin{equation}\label{E-5-1}
\dot k(t)=-\nabla\mathcal E
\end{equation}
with an initial condition in $\mathbb{R}^{\vert V\vert}$. This leads us to the following theorem.
\begin{theorem}\label{T-5-1}
The solution $k(t)$ of \eqref{E-5-1} exists for all time. Moreover, the following two statements are equivalent:
\begin{itemize}
\item[($i$)] The solution $k(t)$ converges as $t\to+\infty$.
\item[($ii$)] The vector $(\hat{L}_i)_{i\in V}$ satisfies \[\hat{L}_i>0,\quad\sum_{i\in I}\hat{L}_i<\pi\vert F_I\vert\]
for each subset $I\subset V$.

\end{itemize}
If either condition $(i)$ or $(ii)$ holds, then the solution $k(t)$ of \eqref{E-5-1} converges exponentially fast to the decorated circle packing with the total geodesic curvature $\hat{L}_{i}$ at $i\in V$.
\end{theorem}

Before proving Theorem \ref{T-5-1}, we provide the following technical tool. Readers may refer to the work of Ge and Xu~\cite[Lemma 4.6]{MR3806817} or Ge, Hua and Zhou~\cite[Lemma. 3.6]{MR4235204} for the proofs.
\begin{lemma}\label{L-5-1}
Suppose $f(x)$ is a $C^{1}$ smooth convex function on $\mathbb{R}^{n}$ with $\nabla f(x_{0})=0$ for some $x_{0}\in\mathbb{R}^{n}$. Suppose $f(x)$ is $C^{2}$ smooth and strictly convex in a neighborhood of $x_{0}$. Then the following statements hold:
\begin{itemize}
\item[($a$)] $\nabla f(x)\neq0$ for any $x\notin\mathbb{R}^{n}\setminus\{x_{0}\}$.
\item[$(b)$] Then $\lim_{\Vert x\Vert\to+\infty}f(x)=+\infty$.
\end{itemize}
\end{lemma}
\begin{proof}[\textbf{Proof of Theorem \ref{T-5-1}}]
The classical ODE theory indicates that the solution exist for $t\in[0,\varsigma)$ and $\Vert k(t)\Vert\to+\infty$ as $t\to\varsigma$ if $\varsigma\neq+\infty$. Meanwhile, we notice that
\[\Vert\dot k(t)\Vert^2=\Vert\nabla\mathcal E\Vert^2=\sum_{v\in V}(L_i-\hat L_i)^2\leq\sum_{v\in V}(C+\vert\hat L_i\vert)^2<+\infty\]
where $C=(\max_{v\in V}\deg v)\pi$ serves as an upper bound for the combinatorial geodesic curvature at each vertex in
$V$, as guaranteed by Lemma \ref{L-4-2}. Then the solution of (\ref{E-5-1}) exists for all time.

Then we prove ($ii$) from ($i$). Suppose $k(t)\to\overline k$ as $t\to\infty$. It follows that $\mathcal E(k(t))\to\mathcal E(\overline k)$. By the mean value theorem, there exists $\xi_n\in(n,n+1)$ such that
\[
\mathcal E(k(n+1))-\mathcal E(k(n))=\frac{d}{dt}\mathcal E(k(\xi_n))\to 0
\]
as $n\to+\infty$. 
In light of $k(\xi_n)\to \overline k$, we know that $\nabla\mathcal E(\overline k)=0$, which implies that $\nabla \Phi(\overline k)=\hat{L}$. From Lemma \ref{L-4-2}, we conclude the prescribed numbers in (\ref{E-5-1}) must satisfy the inequality in ($ii$).

Then, we aim to establish the exponential convergence of the solution from ($ii$).
Let $k(t)$ be the solution of (\ref{E-5-1}).  It is easy to see $\{k(t)\}_{t\geq 0}$ is bounded from above because (\ref{E-5-1}) is a negative gradient flow. Meanwhile, we know $\{k(t)\}_{t\geq 0}$ lies in a compact of $\mathbb{R}^{|V|}$ because $\mathcal E$ is bounded from below and proper, as supported by Lemma \ref{L-5-1}.
Let us define
\[\Upsilon(k(t))=\Vert\nabla\mathcal E(k(t))\Vert^2.\]
We denote the Hessian of $\mathcal E$ by $\mathrm{M}$, it is shown to be positive-definite in the proof of Lemma~\ref{L-4-1}. A direct calculation gives that
\begin{equation}\label{E-5-3}\begin{aligned}
\frac{d\Upsilon}{dt}&=-2\nabla\mathcal E^{\mathrm{T}}\mathrm{M}\nabla\mathcal E,\\
&\leq-2\lambda\Vert\nabla\mathcal E\Vert^2\\
&=-2\lambda \Upsilon,
\end{aligned}\end{equation}
where $\lambda$ is a positive number less than the minimal eigenvalue of $\mathrm{M}$.
The existence of $\lambda$ follows from the fact that $\{k(t)\}_{t\geq 0}$ lies in a compact of $\mathbb{R}^{|V|}$.
As a result of (\ref{E-5-3}), we obtain
\[\Upsilon\left(k(t)\right)\leq\Upsilon(k(0))e^{-2\lambda_0 t},\]
which implies the exponential convergence of the solution.
\end{proof}
\section{Appendix}
This section is devoted to the proof of Lemma \ref{L-3-3}($ii$) and an auxiliary lemma employed in the proof of  Lemma \ref{L-3-3}($i$).
\begin{lemma}\label{L-6-1}
Let $\mathcal{C}$ be a smooth curve in $\mathbb{H}^2$ of constant curvature $1$, and let $\gamma$ be a geodesic in $\mathbb{H}^2$ intersecting $\mathcal{C}$ at two points with intersection angle $\alpha$. Then the arc length of $\mathcal{C}$ between the two points of intersection is equal to $2\tan\alpha$.
\end{lemma}

\begin{proof}
Without loss of generality, we work in the upper half-plane model and assume that
\[\mathcal C=\{(x,y)\in\mathbb{R}^2\ |\ y=1\}.\]
Then $\gamma$ can be parametrized by
\[
\left\{
\begin{aligned}
&x(t)=\tan\alpha(2t-1),\\
&y(t)=1,
\end{aligned}
\right.
\]
for $0\leq t\leq 1$. The hyperbolic length of the arc is then given by
\[l=\int_0^1\frac{\sqrt{x'(t)^2+y'(t)^2}}{y(t)}dt=2\tan\alpha.\]
\end{proof}
\begin{proof}[\textbf{Proof of Lemma \ref{L-3-3}($ii$)}]

The $C^1$-smoothness of $h$ is obvious. To prove the $C^1$-smoothness of $g$,  it suffices to verify its smoothness when $k_1=1$. A direct computation gives that
\begin{equation}\label{E-6-1}\begin{aligned}
g_x(\kappa_1,\kappa_2)=\begin{cases}
\frac{2\kappa_1^2\kappa_2b^2_{\kappa_1}}{\kappa_1^2+\kappa_2^2-1}-2b_{\kappa_1}^3\cot^{-1}b_{\kappa_1}\kappa_2, & \kappa_1>1,\\[1em]
\frac{-2\kappa_1^2\kappa_2b^2_{\kappa_1}}{\kappa_1^2+\kappa_2^2-1}+2b_{\kappa_1}^3\coth^{-1}b_{\kappa_1}\kappa_2, & 0<\kappa_1<1,\\
\end{cases}
\end{aligned}\end{equation}
where $b_{\kappa_s}=\frac{1}{\sqrt{|\kappa_s^2-1|}}$. First, we want to show that
\[\lim_{\kappa_1\to 1^+}g_x(\kappa_1,\kappa_2)=\lim_{\kappa_1\to 1^-}g_x(\kappa_1,\kappa_2).\]
We shall use the Taylor series to compute the limit.
Based on the substitution $\kappa_1=1+\varepsilon$, one can compute that
\begin{equation}\label{E-6-2}\begin{aligned}
b_{\kappa_1}
&=\frac{1}{\sqrt{2\varepsilon}}\Bigl(1 - \tfrac{\varepsilon}{4} + O(\varepsilon^2)\Bigr),\\
\cot^{-1}\bigl(b_{\kappa_1}\kappa_2\bigr)
&=\frac{\sqrt{2\varepsilon}}{\kappa_2}\Bigl(1 + \tfrac{\varepsilon}{4}\Bigr)
- \frac{(2\varepsilon)^{3/2}}{3\,\kappa_2^3}
+ O(\varepsilon^{5/2}).
\end{aligned}\end{equation}
Let us denote
\[\begin{aligned}T_1&=\frac{2\kappa_1^2\kappa_2b^2_{\kappa_1}}{\kappa_1^2+\kappa_2^2-1},\quad T_2=-2b_{\kappa_1}^3\cot^{-1}b_{\kappa_1}\kappa_2.\end{aligned}\]
Then we obtain
\[\begin{aligned}T_1
&=\frac{1}{\varepsilon \kappa_2}-\frac{2}{\kappa_2^3}+\frac{3}{2\kappa_2}+o(\varepsilon),\\
T_2
&=-\frac{1}{\varepsilon \kappa_2}+\frac{2}{3\kappa_2^3}+\frac{1}{2\kappa_2}+o(\varepsilon).\end{aligned}\]
It follows that
\[\lim_{\kappa_1\to 1}g_x(\kappa_1,\kappa_2)=\lim_{\varepsilon\to 0} (T_1+T_2)=\frac{2}{\kappa_2}-\frac{4}{3\kappa^3_2}.\]
Noting that the Taylor expansion of $g_x$ for $0<\kappa_1<1$ is structurally identical to the above case; we henceforth omit its detailed derivation. Then we obtain
\begin{equation}\label{E-6-3}\lim_{\kappa_1\to 1}g_x(\kappa_1,\kappa_2)=\frac{2}{\kappa_2}-\frac{4}{3\kappa^3_2}.\end{equation}
Next, we show that $g$ is differentiable at $\kappa_1=1$. Set $\kappa_1=1+\varepsilon$. By (\ref{E-6-2}), we derive

\[
g(1+\varepsilon, \kappa_2)
= 2\,b_{\kappa_1}\,(1+\varepsilon)\,\cot^{-1}\bigl(b_{\kappa_1}\kappa_2\bigr)
= \frac{2}{\kappa_2}
+ \varepsilon\Bigl(\frac{2}{\kappa_2} - \frac{4}{3\,\kappa_2^3}\Bigr)
+ O(\varepsilon^2),
\]
which yields that
\[
g_+'(1,\kappa_2)
= \lim_{\varepsilon \to 0^+} \frac{g(1+\varepsilon, \kappa_2) - g(1, \kappa_2)}{\varepsilon}
= \frac{2}{\kappa_2} - \frac{4}{3\,\kappa_2^3}.
\]
By a similar derivation, one obtains
\[g_-'(1,\kappa_2)=\frac{2}{\kappa_2} - \frac{4}{3\,\kappa_2^3}.\]
So far, we have shown that
\[g_+'(1,\kappa_2)=g_-'(1,\kappa_2)=\lim_{\kappa_1\to 1}g_x(\kappa_1,\kappa_2),\]
which means $g$ is $C^1$-smooth at $\kappa_1=1$.
\end{proof}

\section{Acknowledgments}
Te Ba is supported by NSF of China (No. 11631010). Guangming Hu is supported by NSF of China (No. 12101275) and Natural Science Research Start-up Foundation of Recruiting Talents of Nanjing University of Posts and Telecommunications (No. NY224040).  Yu Sun is supported by University level natural science foundation of Nanjing Institute of Technology  No.3534113223051. The authors would like to thank Xin Nie, Xu Xu and Ze Zhou for helpful discussions.

\bibliographystyle{siam}

\Addresses
\end{document}